\documentclass[11pt]{amsart}
\usepackage{amsmath}
\usepackage{graphics}      
\theoremstyle{plain}
\newtheorem{theorem}{Theorem}[section]
\newtheorem{lemma}[theorem]{Lemma}

\newtheorem{theo}{Theorem}
\newtheorem{coro}[theo]{Corollary}

\theoremstyle{remark}

\def \CPb {\overline{\mathbf{CP}}^{\,2}}
\def \CP {\mathbf{CP}^2}
\def \RP {\mathbf{RP}^4}

\def \Z {\mathbf{Z}}

\def \KS{\mathrm{KS}}

\begin{document}

\title{A fake smooth $\CP\#\RP$}
\author[Daniel Ruberman]{Daniel Ruberman}
\address{Department of Mathematics, Brandeis University \newline
\hspace*{.375in}Waltham, Massachusetts 02254}
\email{\rm{ruberman@binah.cc.brandeis.edu}}
\thanks{The first author was partially supported NSF Grant DMS9650266 and
the second
author by NSF Grant DMS9626330}
\author[Ronald J. Stern]{Ronald J. Stern}
\address{Department of Mathematics, University of California \newline
\hspace*{.375in}Irvine,  California 92697}
\email{\rm{rstern@math.uci.edu}}
\begin{abstract}
We show that the manifold
$*\RP \# *\CP$, which is homotopy equivalent but not homeomorphic to
$\RP \# \CP$, is in fact smoothable.
\end{abstract} 

\maketitle

\section{Introduction\label{Intro}}
In Kirby's problem list~\cite[Problem 4.82]{kirby:problems96} and in a
recent lecture at MSRI,
P.~Teichner raised the question of the smoothability of a certain non-orientable
$4$-manifold.  In this note we show that the manifold in question, denoted
$*\RP \# *\CP$, which is homotopy equivalent but not homeomorphic to
$\RP \# \CP$, is in fact smoothable. 
The smooth model we construct  will have the additional 
property that its universal cover is
diffeomorphic to
$\CP \#
\CPb$.
 To
describe the manifold in question, we remind the reader that one of the first
consequences of Freedman's simply-connected surgery theory was a construction of
a manifold
$*\CP$, sometimes called
$\mathrm{CH}$ in honor of Chern, which is homotopy equivalent but not homeomorphic to
$\CP$.  The manifold $*\CP$ is not smoothable for classical reasons: it has non-trivial
Kirby-Siebenmann invariant
$\KS\in\Z_2$.  Given any simply-connected {\sl non-spin}
manifold $M$, a similar construction produces a homotopy equivalent
`$*$--partner' $*M$ with opposite Kirby-Siebenmann invariant~\cite{Teich}.   In 1983, the
first author~\cite{ruberman:involution} constructed what is in effect the
$*$-partner of $\RP$.  The connected sum $*\CP \# *\RP$ has trivial
$\KS$-invariant and so might expected to be smoothable; on the other
hand~\cite{hambleton-kreck-teichner:z2} it is not homeomorphic to $\CP \#\RP$.

\begin{theo}\label{smooth}
The manifold $*\CP \# *\RP$ has a smooth structure.   Moreover, it has a smooth
structure such that its universal cover is diffeomorphic to $\CP \# \CPb$
\end{theo}
The classification~\cite{hambleton-kreck-teichner:z2} of non--orientable
manifolds with $\pi_1=\Z_2$ implies that such manifolds which have $b_2 >1$ are
smoothable if and only if their Kirby-Siebenmann invariant vanishes.  Together
with theorem~\ref{smooth} this yields:
\begin{coro}
Let $X$ be a closed non--orientable $4$-manifold with $\pi_1(X) = \Z_2$.  Then
$X$ has a smooth structure if and only if $\KS(X) = 0$.
\end{coro}

\section{Construction of the manifold}

The proof of Theorem~\ref{smooth} is constructive; we will find a smooth
manifold homeomorphic to $*\CP \# *\RP$.  The construction uses a homology
sphere satisfying the conclusion of the following lemma, whose proof will be
given in the next section. 
\begin{lemma}\label{seifert}
There is a homology $3$--sphere $\Sigma^3$ with the following properties.
\begin{itemize}
\item[(i)]  $\Sigma$ is obtained by $\pm 1$ surgery on a knot $K$ in $S^3$.
\item[(ii)]  The Rohlin invariant $\mu(\Sigma) = 1 \pmod{2}$.
\item[(iii)]  $\Sigma$ admits a free, orientation preserving involution $\tau$,
which is isotopic to the identity.
\end{itemize}
\end{lemma}
Different $\Sigma$'s could in principle give rise to different smooth
structures on $*\CP \# *\RP$, but we know of no way to tell them apart.  The
situation is quite analogous to that for the fake $\RP$'s constructed
in~\cite{fs:involution}.
\begin{proof}[Proof of Theorem~\ref{smooth}]  Let $\Sigma$ be a homology
$3$--sphere as described in the lemma; choose an orientation on $\Sigma$ so that
it becomes surgery on a knot with coefficient $= +1$.  Items {\sl (i)} and
{\sl (ii)} are the ingredients in Freedman's
construction~\cite{freedman:simply-connected} of
$*\CP$. That is, let $Y$ be the result of adding a $2$--handle to $B^4$ along
$K$, with framing $1$, then $\partial Y = \Sigma$ and
$$
*\CP = Y \cup_\Sigma \Delta^4
$$ 
where $\Delta^4$ is a contractible
$4$-manifold with boundary $-\Sigma$.  (The sign of the framing
is not really important, for the difference between $*\CP$ and $*\CPb$ will
disappear when we connect sum with $*\RP$.)  The non-trivial $\mu$-invariant is
readily identified with the Kirby-Siebenmann invariant of $*\CP$.

Now items {\sl (ii) } and
{\sl (iii)} are exactly the ingredients for the construction of $*\RP$ given 
in~\cite{ruberman:involution}, i.e.
$$
*\RP = \Delta^4/(x \in \Sigma \sim \tau(x)) = 
(\Sigma/\tau\ \widetilde{\times}\ \mathbf{I})
\cup_\Sigma \Delta^4$$
(The authors of~\cite{hambleton-kreck-teichner:z2} seem to have been unaware of
this earlier construction of  $*\RP$; compare the discussion
in~\cite[Problem 4.74]{kirby:problems96}.)

Let $X$ be the smooth manifold obtained as the union of $Y$ and the mapping
cylinder of the orbit map of the free involution $\tau$ on $\Sigma$, i.e.
$$
X= Y \cup_\Sigma (\Sigma/\tau\ \widetilde{\times}\ \mathbf{I}) = 
Y/(x \in \Sigma \sim \tau(x)).
$$
Then $X$ is manifestly smooth, and we claim that it is homeomorphic to 
$*\CP\#*\RP$.  This seems quite plausible, for the construction amounts to
performing a sort of connected sum, where instead of removing disks and gluing,
we remove the `pseudo-disc' $\Delta^4$ and glue up.  Unfortunately, we do not
know an elementary proof, and must appeal to the homeomorphism classification
theorem of~\cite{hambleton-kreck-teichner:z2}.  

According to that work, the manifold $*\CP \# *\RP$ is distinguished among
non-orientable manifolds with $\pi_1 = \Z_2$ by having $b_2=1$, trivial
Kirby-Siebenmann invariant, and by a codimension-$2$ $\mathrm{Pin}^c$
Arf-invariant.  (The other possible manifolds, up to homeomorphism, 
with the same homology
are
$\CP \#\RP$, $*\CP \# \RP$, and $\CP \# *\RP$.)
The Arf-invariant, whose value for $*\CP \# *\RP$ is $\pm 3 \pmod{8}$,
is that of a surface pulled back from
$\mathbf{CP}^N$ via a map $\varphi:X \to
\mathbf{CP}^{N+1}$ which classifies $c_\Phi$ of the (primitive)
$\mathrm{Pin}^c$ structure $\Phi$.  

A (topological) $\mathrm{Spin}^c$ structure on $*\CP$ also determines such a map,
say $\varphi'$; it is easy to see that (in terms of the decomposition of
$*\CP$ given above) that $\varphi'$ can be taken to be smooth on $Y$, and
constant on $\Delta^4$.  To be more concrete, the dual surface $F$ could be taken
as a Seifert surface of $K$, capped off in the $2$-handle.
The Arf invariant of $F$ (in $*\CP$) is
$4 \pmod 8$, as can be seen from this description of $F$, or by using Rohlin's
theorem as in~\cite{hambleton-kreck-teichner:z2}.

The $\mathrm{Pin}^c$ structure on $*\RP$ has for its characteristic class the
non-trivial class in $H^2(*\RP;\Z)$.  This class is `dual' to a surface in
$*\RP$ which again may be assumed to lie in $*\RP - \Delta^4$.  By the homotopy
invariance of the Arf-invariant for $\mathrm{Pin}^-$  structures, 
$\mathrm{Arf}(*\RP) \equiv  \mathrm{Arf}(\RP) \equiv \pm 1 \pmod 8$.
There is a unique
$\mathrm{Pin}^c$ structure on $\Sigma$, so the  $\mathrm{Pin}^c$ structures on
$*\RP - \Delta^4$ and $*\CP - \Delta^4$ glue up to give a $\mathrm{Pin}^c$
structure $\Phi_X$ on $X$.  The characteristic class $c_{\Phi_X}$ is clearly
dual to the disjoint union of surfaces lying in the two pieces of $X$, so the
Arf invariant is $\pm 4 \pm 1 \equiv \pm 3 \pmod 8$, just as for $*\CP \# *\RP$. 
Since $X$ is smooth, its Kirby-Siebenmann invariant is trivial, and so
$X$ is homeomorphic to $*\CP\# *\RP$.

The additional remark about the universal cover of $X$ being standard may be
seen as follows(cf.~\cite{fs:involution}).  
By the construction of $X$, its cover
$\widetilde{X}
\cong Y \cup_\tau \overline{Y} \cong Y \cup \overline{Y}$ since $\tau$ is
isotopic to the identity.  On the other hand, $Y \cup \overline{Y}$ is obtained
by adding two $2$-handles to $B^4$, together with a $4$-handle.  The first is
added along
$K$, with framing
$1$, and the second is added along a meridian of $K$, with framing $0$.  (This
is a standard argument in handle theory, see for
example~\cite{kirby:4-manifolds}.)  It is then easy to unknot $K$, by
repeatedly sliding over the $0$-framed handle, resulting in a standard
picture of $\CP\# \CPb$.
\end{proof}

\section{Proof of Lemma~\ref{seifert}}

In this section, we give two examples of homology spheres satisfying the
conclusions of Lemma~\ref{seifert}.  Both examples are Brieskorn
spheres, i.e. Seifert-fibered homology spheres of the form $\Sigma(p,q,r)$,
where $p,q,$ and
$r$ are relatively prime odd numbers.   The involution $\tau$ is nothing more
than multiplication by  $-1 \in S^1$ in the natural circle action on
$\Sigma(p,q,r)$.  The condition that the numbers $p,q,$ and $r$ be odd
guarantees that $\tau$ is free; since $-1$ is contained in a circle, the
involution is isotopic to the identity.

There are many Brieskorn spheres which are integral surgery on a
knot--for some examples see~\cite{kalliongis:seifert,miyazaki:seifert} or
adapt the technique of~\cite{casson-harer}.  For most of these constructions
one of the
indices turns out to be even.  One  construction is given below, where
it is shown that adding a handle (along the curve denoted $\gamma$)
to the Brieskorn sphere $\Sigma(5,9,13)$ yields  $S^3$.  Turning the picture
upside down shows that $\Sigma(5,9,13)$ is integral surgery on a knot in $S^3$. 
As remarked in the proof of Theorem~\ref{smooth}, it doesn't matter whether the
coefficient is positive or negative.  Again, the $\mu$-invariant is $1$ (from
the picture just after blowing down the first $-1$ curve), so this example
proves the lemma.


Another construction from the
literature which provides Seifert fibered spaces is $rs(p+q)^2 + pq$
surgery on the knot denoted
$K_{p,q}(r,s)$ in the recent paper~\cite[\S 9]{miyazaki:seifert}.  
Choosing $p = -13$, $q=23$, $r= 3$, and $s=1$ gives the homology sphere
$\Sigma(3,13,23)$ as $ +1$ surgery on a hyperbolic
knot.  Since $\mu(\Sigma(3,13,23)) =1$, this manifold gives an example which
yields the proof of Lemma~\ref{seifert}.   This is the only example of a
$\mu$-invariant $1$ homology sphere constructible by this method found by a
moderately long computer search.
It is possible to give a Kirby-calculus proof that $\Sigma(3,13,23)$ is
surgery on a knot similar to the one for $\Sigma(5,9,13)$; aficionados of the
subject may wish to check if the knot is the same as the one in the knot from the
paper~\cite{miyazaki:seifert}.

\begin{center}
\leavevmode  
\includegraphics{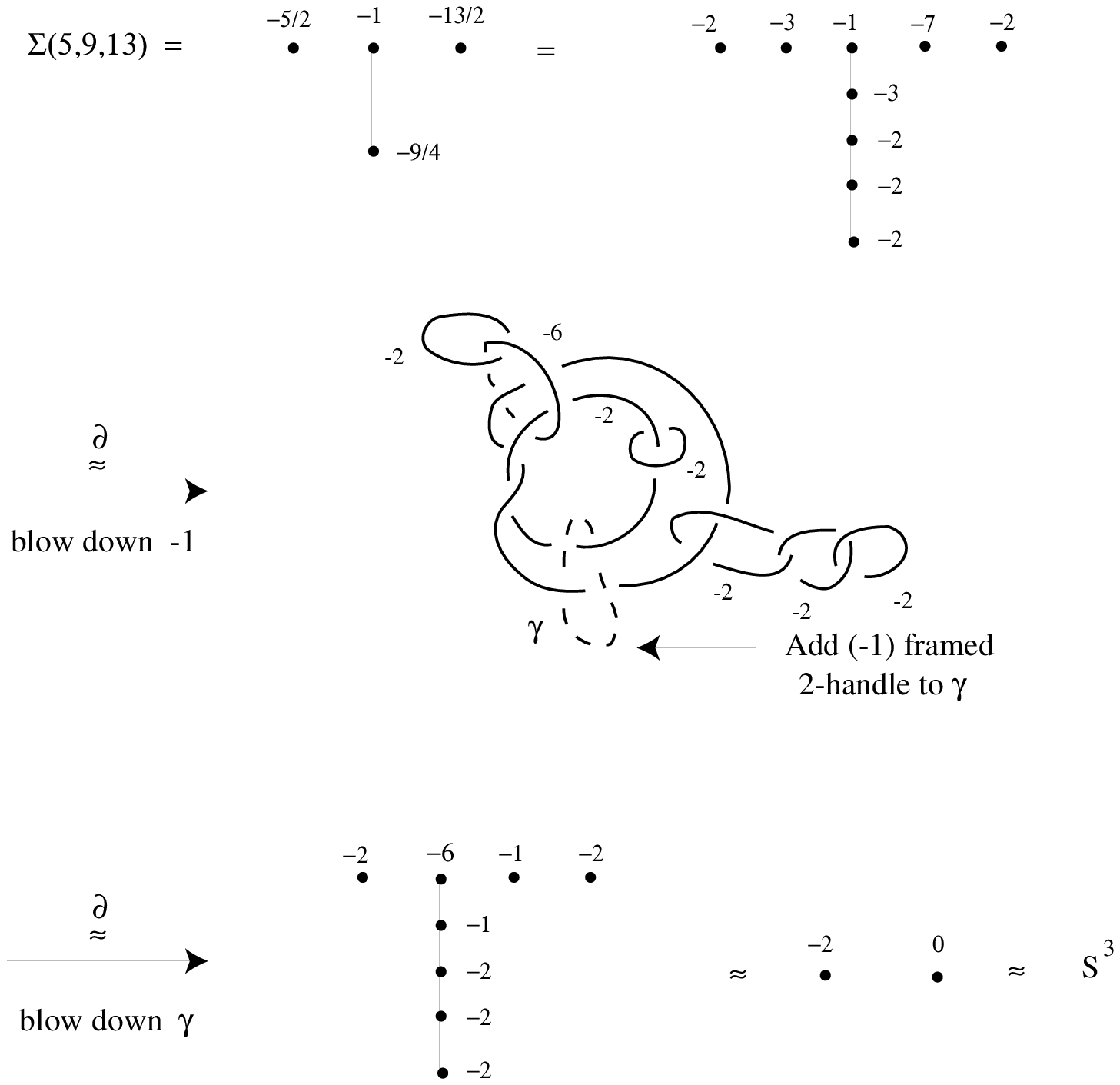}
\end{center}



\begin{thebibliography}{HKT94}

\bibitem[CH81]{casson-harer}
A.~Casson and J.~Harer, \emph{Some homology lens spaces which bound rational
  homology balls}, Pacific J.\ Math. \textbf{96} (1981), 23--36.

\bibitem[Fre82]{freedman:simply-connected}
M.~Freedman, \emph{The topology of four--dimensional manifolds}, J. Diff.\ Geo.
  \textbf{17} (1982), 357--432.

\bibitem[FS81]{fs:involution}
R.~Fintushel and R.~Stern, \emph{An exotic free involution on {$S\sp{4}$}},
  Annals of Math. \textbf{113} (1981), 357--365.

\bibitem[HKT94]{hambleton-kreck-teichner:z2}
I.~Hambleton, M.~Kreck, and P.~Teichner, \emph{Nonorientable {$4$}-manifolds
  with fundamental group of order {$2$}.}, Trans.\ A.M.S. \textbf{344} (1994),
  649--665.

\bibitem[Kir89]{kirby:4-manifolds}
R.C. Kirby, \emph{Topology of $4$-manifolds}, Lecture Notes in Math., vol.
  1374, Springer-Verlag, 1989.

\bibitem[Kir97]{kirby:problems96}
R.C. Kirby, \emph{Problems in low--dimensional topology}, Geometric Topology
  (W.~Kazez, ed.), American Math.\ Soc./International Press,
  Providence, 1997.

\bibitem[KT90]{kalliongis:seifert}
John Kalliongis and Chichen~M. Tsau, \emph{Seifert fibered surgery manifolds of
  composite knots}, Proc. Amer. Math. Soc. \textbf{108} (1990), 1047--1053.

\bibitem[MM97]{miyazaki:seifert}
Katura Miyazaki and Kimihiko Motegi, \emph{Seifert fibred manifolds and {D}ehn
  surgery}, Topology \textbf{36} (1997), 579--603.

\bibitem[Rub84]{ruberman:involution}
D.~Ruberman, \emph{Equivariant knots of free involutions of ${S}^4$}, Top.\
  Appl. \textbf{18} (1984), 217--224.

\bibitem[Teich96]{Teich}
P.~Teichner, \emph{On the star-construction for topological 4-manifolds}, 
Geometric Topology (W.~Kazez, ed.), American Math.\ Soc./International Press,
  Providence, 1997, 300--312.
\end{thebibliography}
\providecommand{\bysame}{\leavevmode\hbox to3em{\hrulefill}\thinspace}

\end{document}